\documentclass[11pt,a4paper]{article}
\usepackage[T1]{fontenc}
\usepackage{lmodern}
\usepackage{microtype}
\usepackage[margin=2.45cm]{geometry}
\usepackage{amsmath,amssymb,amsthm,mathtools,mathrsfs}
\usepackage{authblk}
\usepackage[colorlinks=true,linkcolor=blue,citecolor=blue,urlcolor=blue]{hyperref}
\hypersetup{
 pdftitle={Local Limit Theorem for Branching Processes in Varying Environments},
 pdfauthor={Zhiyang Xia and Mei Zhang},
 pdfsubject={Uniform local limits for branching processes in varying environments}
}
\numberwithin{equation}{section}
\allowdisplaybreaks
\theoremstyle{plain}
\newtheorem{theorem}{Theorem}[section]
\newtheorem{lemma}[theorem]{Lemma}
\newtheorem{proposition}[theorem]{Proposition}

\theoremstyle{remark}

\newcommand{\E}{\mathbb E}
\newcommand{\Prob}{\mathbb P}
\newcommand{\Var}{\operatorname{Var}}
\newcommand{\Nzero}{\mathbb N_0}
\newcommand{\norm}[1]{\left\lVert#1\right\rVert_{1}}
\newcommand{\dd}{\mathbb D}
\newcommand{\HH}{\mathcal H}
\newcommand{\RR}{\mathcal R}
\newcommand{\condref}[2]{\hyperlink{cond:#1}{#2}}
\begin{document}

\title{Local Limit Theorem for Branching Processes\\
in Varying Environments
}

\author[a]{Zhiyang Xia\thanks{Email: 202521130107@mail.bnu.edu.cn}}
\author[b]{Mei Zhang\thanks{Corresponding author: meizhang@bnu.edu.cn}}
\affil[a,b]{Laboratory of Mathematics and Complex Systems (Ministry of Education), School of Mathematical Sciences, Beijing Normal University\\
Beijing 100875, P. R. China}
\date{}
\maketitle

\begin{abstract}
Let $Z_{i,n}$ be the population at time $n$ descended from one particle
at time $i$ in a branching process in a varying environment.
Under two-sided bounds on products of
offspring means, an integrable uniform second-moment tail condition, and
nondegeneracy and unit positive-support gcd of all weak limits, we prove
$\Prob(Z_{i,n}=k)\le C(n-i+1)^{-2}$ uniformly in $0\le i<n$ and $k\ge1$,
and obtain a local limit theorem as $n-i,k\to\infty$ with
$k\le L_0(n-i)$, without third-moment assumptions. These results extend the classical local probability estimates for critical Galton--Watson processes to varying environments.
\end{abstract}

\bigskip


\section{Introduction and main results}

Let $\mathbf{Z^0}=\{Z_n^0,n\ge0\}$ be a Galton--Watson (GW)
branching process defined by
\begin{align}\label{Def Zn0}
Z_n^0=\sum_{i=1}^{Z_{n-1}^0}\zeta_{ni},\qquad n\ge1,\qquad Z_0^0=1,
\end{align}
where $\{\zeta_{ni}:n,i\ge1\}$ is a sequence of independent and
identically distributed (i.i.d.) random variables with a nondegenerate
common distribution $\{p_j:j\ge0\}$. Denote
$m:=\sum_{j=0}^\infty jp_j$ and
$0<B:=\sum_{j=2}^\infty j(j-1)p_j<\infty$. The process is called
supercritical, critical, or subcritical according as $m>1$, $m=1$, or
$m<1$. Its lattice span is
\[
 d:=\gcd\{i-j:i\ne j,\ p_ip_j>0\},
\]
and the case $d=1$ is called aperiodic.

A classical local limit theorem for the critical Galton--Watson
process $\mathbf{Z^0}$ was established by Kesten, Ney and Spitzer
\cite{KNS}. In the aperiodic case, under the additional assumption
$\sum_{j=2}^\infty p_j j^2\log j<\infty$, they proved that, as
$n\to\infty$ and $k\to\infty$ with $k/n$ bounded,
\begin{equation}\label{eq:KNS}
 \Prob(Z_n^0=k)=\frac{4}{B^2n^2}
       \exp\!\left(-\frac{2k}{Bn}\right)(1+o(1)).
\end{equation}
For $d>1$, the leading term is multiplied by $d$ on admissible lattice
points. Nagaev and Vakhtel \cite{NagaevVakhtel} subsequently
established \eqref{eq:KNS} under the finite-variance assumption alone.
For further background on local limit theorems, including those for
critical Bellman--Harris processes, see \cite{NagaevVakhtel} and the
references therein. Limit theorems, including local probability
estimates, have also been studied extensively for critical GW processes
with immigration; see, for example,
\cite{P69,P72,LZ2021,LZ2024}. Related statistical inference for GW
processes is treated in \cite{P75}.

A branching process in a varying environment (BPVE) allows the
offspring law to change with generation. Such processes may have a
finite positive limiting population with positive probability
(Lindvall \cite{Lindvall1974}). Early work on their classification
into subcritical, critical, and supercritical regimes includes Jagers
\cite{Jagers1974}. MacPhee and Schuh \cite{MacPheeSchuh} studied
examples with essentially constant offspring means but different growth
rates. More recently, Kersting \cite{Kersting2020}
developed a classification of BPVEs under a uniform relative moment
condition. Under a stronger relative tail condition, distinct from our
mean-product assumption \condref{B}{(B)}, he also established
Kolmogorov-type survival asymptotics and a Yaglom-type conditional
exponential limit theorem in the critical regime
\cite[Theorem~4]{Kersting2020}. The corresponding scale is the
conditional expectation of the population given survival, which plays
the role of $Bn/2$ in \eqref{eq:KNS}. For critical BPVEs satisfying
the same relative tail condition, the analogy with the critical
Galton--Watson process was extended to the genealogical level in
\cite[Theorems~2 and~3]{KerstingGeneal}. After an
environment-dependent rescaling of time, the most recent common
ancestor given survival and the corresponding reduced process exhibit
the classical limiting behaviour. Mitov and Omey
\cite{MitovOmey2014} studied critical processes with varying offspring
variance and immigration. Later, Gonz\'alez, Kersting, Minuesa and del
Puerto \cite[Theorem~4]{Gonzalez2019} obtained gamma limits under
suitable normalization for processes with generation-dependent
immigration. The BPVE results recalled above concern survival
probabilities, weak limits, and genealogical structure rather than
local probabilities.

In this paper, we study local probabilities for a branching process
in a varying environment $\mathbf Z$, where $Z_{i,n}$ is the population
at time $n$ descended from a single particle at time $i$. We consider
environments in which the products of offspring means are bounded above
and away from zero. We establish an upper bound for
$p_{i,n}(k)=\Prob(Z_{i,n}=k)$ uniformly in the starting time $i$, and
obtain a corresponding local limit theorem as $h=n-i\to\infty$ and
$k\to\infty$ with $k\le L_0h$ for fixed $L_0>0$; see
Theorem~\ref{thm:main}.

In a constant critical environment, $m_r=1$ and $f_r=f$ for every $r$,
so $M_{0,n}=1$ and $a_{0,n}=Bn/2$; hence the asymptotic in
Theorem~\ref{thm:main} is exactly \eqref{eq:KNS}. Nagaev and Vakhtel
\cite{NagaevVakhtel} proved the fixed-environment result under finite
variance alone and also established
$\sup_{n,k\ge1}n^2\Prob(Z_n^0=k)<\infty$. The new feature in the
varying-environment setting is that both the $O((n-i+1)^{-2})$ bound
and the local asymptotic are uniform over all starting times $i$ and
elapsed times $n-i$, with $M_{i,n}$ and $a_{i,n}$ adapting the
normalization to the offspring laws on each interval.

For $Z_0=1$, our assumptions imply almost sure extinction and a
conditional mean of order $n$:
$\E[Z_n\mid Z_n>0]\sim a_{0,n}\asymp n$. These conclusions follow from
Lemmas~\ref{lem:moments}, \ref{lem:survival}, and
\ref{lem:linear-fractional}. One difference from Kersting's framework
is that neither the statement nor the proof of our main result requires
his regularity condition \condref{A}{(A)}. In Kersting~\cite{Kersting2020}, this
condition helps to keep the upper and lower shape-function bounds
together uniformly across generations. Our argument
proceeds differently: the survival upper bound combines the truncated
shape-function lower bound in Lemma~\ref{lem:kersting2}, whose
truncation error is controlled by \condref{T}{(C)}, with the uniform
growth of $a_{i,n}$. In this way, the local conclusions are obtained
without imposing \condref{A}{(A)}.  If condition
\condref{A}{(A)} is additionally imposed, the process is regular in
Kersting's sense, and \cite[Proposition~1]{Kersting2020} places it in
his critical class.

We now define the process and introduce our notation. Let
\begin{equation}\label{Def Zn}
 Z_n=\sum_{j=1}^{Z_{n-1}}\xi_{nj},\qquad n\ge1,\qquad Z_0=1,
\end{equation}
where $\{\xi_{nj}:n,j\ge1\}$ are independent non-negative integer-valued
random variables and, for each $n$, the variables
$\{\xi_{nj}:j\ge1\}$ have the same distribution, with probability
generating function $f_n$. Set $\sum_{j=1}^{0}=0$.
Let $Y_r$ have the offspring generating function $f_r$, and write
$m_r=\E Y_r\in(0,\infty)$. Define
\[
 \mu_0=1,\qquad \mu_n=\prod_{r=1}^n m_r,\qquad
 \nu_r=\frac{f_r''(1)}{m_r^2},\qquad
 \beta_r=\frac{f_r''(1)}{m_r}.
\]
For $0\le i\le n$, let
\[
 F_{i,n}=f_{i+1}\circ\cdots\circ f_n,\qquad F_{n,n}(s)=s.
\]
Let $Z_{i,n}$ denote the number of particles at time $n$,
starting from one particle at time $i$. Thus $F_{i,n}$ is the
probability generating function of $Z_{i,n}$. Write
\begin{align}
 & p_{i,n}(k)=\Prob(Z_{i,n}=k),\quad
 q_{i,n}=F_{i,n}(0),\quad Q_{i,n}=1-q_{i,n},\quad
 M_{i,n}=\frac{\mu_n}{\mu_i}.
\label{Qin}\\
& a_{i,n}=\frac{\mu_n}{2}\sum_{r=i+1}^{n}\frac{\nu_r}{\mu_{r-1}},
 \qquad a_{n,n}=0.\label{eq:a-def}
\end{align}

We write $\E[X;A]=\E[X\mathbf1_A]$ for an event $A$.
For a nonnegative integer-valued random variable $Y$ with
$\Prob(Y\ge1)>0$, write
\[
 d_+(Y):=\gcd\{k\ge1:\Prob(Y=k)>0\}.
\]

For comparison with Kersting's regularity condition, we record the
following condition, which is not assumed in our main result:

\hypertarget{cond:A}{(A)} there is a constant $C_A\ge1$ such that, for every $r\ge1$,
\begin{equation}\label{eq:relative-moment-condition}
 \E[Y_r^2;Y_r\ge2]
 \le C_A\E[Y_r;Y_r\ge2]\E[Y_r\mid Y_r\ge1].
\end{equation}

We use the following three conditions. They are in force throughout
Sections~2 and~3.
\begin{enumerate}
\item[\hypertarget{cond:B}{(B)}] The products of the means satisfy
\[
 0<\underline\mu\le\mu_n\le\overline\mu<\infty,
 \qquad n\ge0.
\]
\item[\hypertarget{cond:T}{(C)}] The uniform second-moment tail envelope satisfies
\begin{equation}\label{eq:tail-integral-condition}
 \begin{gathered}
 V(x):=\sup_{r\ge1}\E[Y_r^2;Y_r>x],\qquad x\ge0,\\
 K_T:=\int_0^\infty\frac{V(x)}{1+x}\,\mathrm dx<\infty .
 \end{gathered}
\end{equation}
\item[\hypertarget{cond:D}{(D)}] Every weak limit $Y_*$ of a subsequence
$Y_{r_j}$ with $r_j\uparrow\infty$ satisfies both
\[
 \Var(Y_*)>0,\qquad d_+(Y_*)=1.
\]
\end{enumerate}

The following is the main theorem of the present paper.
\begin{theorem}\label{thm:main}
Under assumptions \condref{B}{(B)}, \condref{T}{(C)},
and \condref{D}{(D)}, there exists a positive constant $C$ such that
\begin{equation}\label{eq:main-upper}
 \sup_{k\ge1}p_{i,n}(k)\le\frac{C}{(n-i+1)^2},
 \qquad 0\le i<n.
\end{equation}
Moreover, for any fixed $L_0>0$,
\begin{equation}\label{eq:main-local}
\lim_{N\to\infty}
 \sup_{\substack{0\le i<n,\ n-i\ge N\\
                   k\in\mathbb N,\ N\le k\le L_0(n-i)}}
 \left|
 \frac{a_{i,n}^{2}}{M_{i,n}}
       \exp\!\left(\frac{k}{a_{i,n}}\right)p_{i,n}(k)-1
 \right|=0.
\end{equation}
The supremum is taken for sufficiently large $N$, for which all scales
$a_{i,n}$ in its range are positive. In particular, for the process
starting from one particle at time $0$, write $Z_n=Z_{0,n}$ and
$a_n=a_{0,n}$. Then, as $n\to\infty$ and $k\to\infty$ with $k/n$
bounded,
\[
 \Prob(Z_n=k)=\frac{\mu_n}{a_n^2}
       \exp\!\left(-\frac{k}{a_n}\right)(1+o(1)).
\]
\end{theorem}

Our proof of Theorem~\ref{thm:main} extends the generating-function
approach of \cite{KNS} to varying offspring laws. The main difficulty is that the offspring law
changes from one generation to the next. We therefore need estimates
for compositions of different generating functions that hold uniformly
over all time intervals. The argument has three main steps.

First, we determine the relevant time scale and estimate the survival
probabilities. Lemma~\ref{lem:survival} shows that $a_{i,n}$ grows
linearly with the length of the time interval and gives upper and lower
survival bounds of order $(n-i+1)^{-1}$. The lower survival bound
follows from moment estimates. For the upper bound, we combine the
truncated shape-function bound in Lemma~\ref{lem:kersting2} with the
scale growth $a_{i,n}$, and condition \condref{T}{(C)} allows us to control the
truncation error uniformly.

Second, we approximate the composed generating functions and control
the remaining correction terms. Starting from the reciprocal
composition identity \eqref{eq:reciprocal}, we approximate $1-F_{i,n}$
by a linear fractional function and estimate the error in
Lemma~\ref{lem:linear-fractional}. Lemmas~\ref{lem:shape-modulus} and
\ref{lem:dini} provide the uniform error bounds needed for this
approximation. A second-order expansion of the one-step factor
(Lemma~\ref{lem:one-step}) then leads to an exact product recursion for
the normalized correction factors in \eqref{eq:correction-def}
(Lemma~\ref{lem:product}). Iterating this recursion gives a product
representation, and Proposition~\ref{prop:correction-compact} supplies
the compactness and uniform coefficient-tail bounds needed below.

Finally, Section~3 expands the remaining geometric factor and extracts
the coefficients of the generating function. Together with the uniform
bounds from the previous steps, this yields both the upper estimate and
the local asymptotic in Theorem~\ref{thm:main}.

The argument also uses the shape-function inequalities and the moment
and reciprocal composition identities of
\cite[Lemmas 1, 4 and 5]{Kersting2020}. For background on generating
functions and convolution methods, see \cite{SpitzerRW,FellerII};
related survival asymptotics in a random environment are studied in
\cite{Kozlov1976}.

The paper is organized as follows. Section~2 collects the auxiliary
lemmas: Subsection~2.1 introduces the coefficient norm on power series
and the second-moment error modulus; Subsection~2.2 establishes uniform
estimates for moments, cumulative scales, and survival probabilities;
Subsection~2.3 establishes the gcd criterion for inversion and studies
the single-step shape functions and their composition errors;
Subsection~2.4 derives the second-order Taylor expansion and the exact
product recursion; and Subsection~2.5 establishes relative compactness
of the family of normalized correction factors. In Section~3 we give the proof of Theorem~\ref{thm:main}.

Throughout the paper, $\Nzero$ is the set of non-negative integers.
The letters $C$ and $c$ denote positive finite constants whose values
may change from line to line. They may depend on the fixed environment,
including the eventual nondegeneracy threshold and the compactness data
of its offspring laws, but never on $i$, $n$, or $k$. A subscript
indicates an additional allowed dependence. In particular, we do not claim that the constants are uniform across all environments with the same numerical bounds in
\condref{B}{(B)} and \condref{T}{(C)}.

\section{Auxiliary lemmas}

This section proves the auxiliary estimates used in Theorem~\ref{thm:main}.

\subsection{The coefficient norm and the second-moment error modulus}

Let $\dd=\{s\in\mathbb C:|s|<1\}$. For a formal power series
$g(s)=\sum_{j\ge0}g_js^j$, define its coefficient norm by
\[
 \norm{g}:=\sum_{j\ge0}|g_j|\le\infty.
\]
All norm and convergence statements below refer to this $\ell^1$
norm of the coefficient sequence. If $\norm g<\infty$, the series
converges absolutely and uniformly on the closed disc
$\overline{\dd}$; in particular, $g$ is analytic on $\dd$ and
continuous on $\overline{\dd}$, and it is uniquely determined by its
values on $(0,1)$. Products of series with finite coefficient norm
again have finite coefficient norm, and absolute convergence gives
the convolution inequality
\begin{equation}\label{eq:convolution}
 \norm{gh}\le\norm g\norm h,
\end{equation}
together with pointwise multiplicativity, $(gh)(s)=g(s)h(s)$ for
every $s\in\overline{\dd}$. A probability generating function $F$
satisfies $\norm F=1$. If $\norm g<\infty$, the composition
$g\circ F$ also has finite coefficient norm, and
\begin{equation}\label{eq:composition-contraction}
 \norm{g\circ F}\le\norm g,\qquad
 \norm{F^j-1}=2\bigl(1-F(0)^j\bigr)\quad(j\ge1).
\end{equation}
Composition also respects products and units:
\begin{equation}\label{eq:composition-product}
 (gh)\circ F=(g\circ F)(h\circ F),\qquad 1\circ F=1,
\end{equation}
since both sides have absolutely summable coefficients and agree
pointwise on $\overline{\dd}$.

The following version of Wiener's inversion theorem originates in
\cite{Wiener1932}; the precise form used here follows from
\cite[Lemma~11.6]{RudinFA} or \cite[Chapter~5]{ZhangGuo}.

\begin{lemma}\label{lem:wiener}
If $g(s)=\sum_{j\ge0}g_js^j$ satisfies $\norm g<\infty$ and has no
zero on the closed unit disc $\overline{\dd}$, then its reciprocal
admits a power series expansion $1/g(s)=\sum_{j\ge0}b_js^j$ with
$\sum_{j\ge0}|b_j|<\infty$.
\end{lemma}

The next lemma gives the uniform tail estimates for $Y_r$ used
throughout the paper.

\begin{lemma}\label{lem:dini}
Under \condref{T}{(C)}, let $K_2:=\sup_{r\ge1}\E Y_r^2.$
Then
$
 K_2\le\frac{K_T}{\log2}<\infty,
$
and
\begin{equation}\label{secondm}
 \sup_{r\ge1}\E[Y_r^2;Y_r>K]=V(K)\longrightarrow0,
                    \qquad K\to\infty.
\end{equation}
Define
\begin{equation}\label{Omegat0}
 \omega(t):=\sup_{r\ge1}\E[Y_r^2\min\{tY_r,1\}],
 \qquad 0\le t\le1.
\end{equation}
This function $\omega(\cdot)$ is nondecreasing, is bounded by $K_2$,
and satisfies $\omega(t)\to0$ as $t\downarrow0$. More precisely,
\begin{equation}\label{eq:omega-bound}
 \omega(t)\le t\int_0^{1/t}V(x)\,\mathrm dx,\qquad
 \int_0^1\frac{\omega(t)}t\,\mathrm dt\le2K_T,
 \qquad 0<t\le1.
\end{equation}
For every fixed $C_0\ge1$,
\begin{equation}\label{eq:dini-sum}
 \sum_{h=0}^\infty\frac1{h+1}
   \omega\!\left(\min\left\{1,\frac{C_0}{h+1}\right\}\right)
 \le K_2(1+\log C_0)+2K_T<\infty.
\end{equation}
In particular, the tails of this series tend to zero. Conversely,
if $K_2<\infty$ and $\int_0^1\omega(t)\,\mathrm dt/t<\infty$, then
\condref{T}{(C)} holds.
\end{lemma}

\begin{proof}
Since the $Y_r$ are non-negative integer-valued, $V(x)=K_2$ for
$0\le x<1$. Consequently,
\[
 K_T\ge\int_0^1\frac{V(x)}{1+x}\,\mathrm dx=K_2\log2.
\]
This proves the bound for $K_2$. Clearly, $V$ is nonincreasing.
If $V(K)\to c>0$ as $K\to\infty$, then $K_T=\infty$, contradicting
\condref{T}{(C)}. Hence $V(K)\to0$.

The uniform integrability just obtained also controls the modulus
$\omega$, which is nondecreasing and bounded by $K_2$.
Splitting at an integer $K\ge1$ gives, for $0<t\le1$,
\begin{equation}\label{eq:omega-prefinal}
 \omega(t)\le tK K_2+V(K).
\end{equation}
Letting first $t\downarrow0$ and then $K\to\infty$ gives
$\omega(t)\to0$.

For every nonnegative real $y$,
\[
 \min\{ty,1\}
 =t\int_0^{1/t}\mathbf1_{\{y>x\}}\,\mathrm dx.
\]
By Tonelli's theorem and taking the supremum over $r$,
\[
 \begin{aligned}
 \omega(t)
 =\sup_r t\int_0^{1/t}\E[Y_r^2;Y_r>x]\,\mathrm dx
 \le t\int_0^{1/t}V(x)\,\mathrm dx .
 \end{aligned}
\]
Integrating this inequality with weight $1/t$ and applying Tonelli's
theorem again gives
\[
 \begin{aligned}
 \int_0^1\frac{\omega(t)}t\,\mathrm dt
 &\le\int_0^1\int_0^{1/t}V(x)\,\mathrm dx\,\mathrm dt\\
 &=\int_0^\infty V(x)\min\{1,x^{-1}\}\,\mathrm dx\\
 &\le2\int_0^\infty\frac{V(x)}{1+x}\,\mathrm dx=2K_T,
 \end{aligned}
\]
and hence \eqref{eq:omega-bound} holds.

To deduce \eqref{eq:dini-sum}, fix $C_0\ge1$ and put
$\psi(t)=\omega(\min\{1,C_0t\})$ for $0<t\le1$.
It is nondecreasing, and a change of variables gives
\[
 \int_0^1\frac{\psi(t)}t\,\mathrm dt
 =\int_0^1\frac{\omega(u)}u\,\mathrm du+\omega(1)\log C_0.
\]
For integers $h\ge1$, monotonicity implies
\[
 \begin{aligned}
 \int_{1/(h+1)}^{1/h}\frac{\psi(t)}t\,\mathrm dt
 \ge\psi(1/(h+1))\log\frac{h+1}{h}
 \ge\frac{\psi(1/(h+1))}{h+1}.
 \end{aligned}
\]
The last inequality uses $\log(1+1/h)\ge1/(h+1)$. Summing and adding
the $h=0$ term $\omega(1)\le K_2$ proves \eqref{eq:dini-sum}.

For the converse, suppose that $K_2<\infty$ and
$\int_0^1\omega(t)\,\mathrm dt/t<\infty$. For $x\ge1$, on $Y_r>x$ one
has $\min\{Y_r/x,1\}=1$, so $V(x)\le\omega(1/x)$. Consequently,
\[
 \begin{aligned}
 \int_0^\infty\frac{V(x)}{1+x}\,\mathrm dx
 \le K_2\log2+
       \int_1^\infty\frac{\omega(1/x)}x\,\mathrm dx
 =K_2\log2+\int_0^1\frac{\omega(t)}t\,\mathrm dt<\infty.
 \end{aligned}
\]
This proves the converse.
\end{proof}

\subsection{Uniform moment estimates and survival probabilities}

For a probability generating function $f$ with $0<m:=f'(1)<\infty$
and $f''(1)<\infty$, put $\nu:=f''(1)/m^2$ and define the
single-step shape function by
\begin{equation}\label{varphi00}
 \varphi(s)=\frac1{1-f(s)}-\frac1{m(1-s)},\quad 0\le s<1,
 \qquad \varphi(1):=\frac{\nu}{2}.
\end{equation}
The following bound is from \cite{Kersting2020}.

\begin{lemma}[{\cite[Lemma~1]{Kersting2020}}]\label{lem:kersting1}
Let $f$ be a probability generating function with
$0<m=f'(1)<\infty$ and $f''(1)<\infty$. Then, for $0\le s\le1$,
\begin{equation}\label{e2.20}
 \frac12\varphi(0)\le\varphi(s)\le2\varphi(1).
\end{equation}
\end{lemma}

Let $\varphi_r$ denote the shape function of $f_r$. The following
truncated lower bound is used in the survival estimate below.

\begin{lemma}[Truncated shape-function lower bound]\label{lem:kersting2}
Under \condref{B}{(B)} and \condref{T}{(C)}, put
$m_*=\underline\mu/\overline\mu$. For every integer $K\ge2$,
$r\ge1$, and $0\le s\le1$,
\begin{equation}\label{eq:relative-shape-comparison}
 \varphi_r(s)\ge\frac{\nu_r}{2K}
                   -\frac{V(K)}{2K m_*^2}.
\end{equation}
\end{lemma}

\begin{proof}
Condition \condref{B}{(B)} gives $m_r\ge m_*>0$, and
\condref{T}{(C)} gives finite second moments. Since
$Y_r(Y_r-1)\le K(Y_r-1)$ on $\{2\le Y_r\le K\}$,
\[
 f_r''(1)=\E[Y_r(Y_r-1)]
 \le K\E[Y_r-1;Y_r\ge1]+V(K).
\]
Moreover, $\Prob(Y_r\ge1)>0$ and $\Prob(Y_r\ge1)\le m_r$, so
\[
 \varphi_r(0)
 =\frac{m_r-\Prob(Y_r\ge1)}{m_r\Prob(Y_r\ge1)}
 =\frac{\E[Y_r-1;Y_r\ge1]}{m_r\Prob(Y_r\ge1)}
 \ge\frac{\E[Y_r-1;Y_r\ge1]}{m_r^2}.
\]
Using Lemma~\ref{lem:kersting1}, we obtain
\[
 \varphi_r(s)\ge\frac12\varphi_r(0)
 \ge\frac{f_r''(1)-V(K)}{2K m_r^2}
 \ge\frac{\nu_r}{2K}-\frac{V(K)}{2K m_*^2},
\]
which proves the claim. The right-hand side is allowed to be negative.
The survival proof will control its accumulated error with a single
fixed choice of $K$.
\end{proof}

\begin{lemma}[Composition formulas for the first two factorial moments]
\label{lem:factorial-moments}
For offspring laws with positive finite means and finite second
moments, and for all $0\le i\le n$,
\begin{equation}\label{eq:process-moments}
 \E Z_{i,n}=M_{i,n},\qquad
 \E[Z_{i,n}(Z_{i,n}-1)]=2a_{i,n}M_{i,n}.
\end{equation}
In particular, $\E Z_{i,n}^2=M_{i,n}(1+2a_{i,n})$.
\end{lemma}

\begin{proof}
These are the composition formulas in
\cite[Lemma~4]{Kersting2020}, applied to the environment beginning
at time $i$. We give the calculation to specify the normalization.
The chain rule yields
\[
 F_{i,n}'(1)=\prod_{r=i+1}^n m_r=M_{i,n}.
\]
Put $A_{i,n}=F_{i,n}''(1)/F_{i,n}'(1)$. Since
$F_{i,n}=f_{i+1}\circ F_{i+1,n}$,
\[
 A_{i,n}=A_{i+1,n}
       +\frac{f_{i+1}''(1)}{m_{i+1}}M_{i+1,n}.
\]
Iterating from $A_{n,n}=0$ gives
\begin{align*}
 A_{i,n}
 &=\sum_{r=i+1}^n\frac{f_r''(1)}{m_r}\frac{\mu_n}{\mu_r}\\
 &=\mu_n\sum_{r=i+1}^n\frac{f_r''(1)}{m_r^2\mu_{r-1}}
   =2a_{i,n}.
\end{align*}
The derivatives of a generating function at one give its first and
second factorial moments. Therefore
$F_{i,n}''(1)=2a_{i,n}M_{i,n}$, proving the two identities.
For $i=n$, they read $\E Z_{n,n}=1$ and
$\E[Z_{n,n}(Z_{n,n}-1)]=0$, as required.
\end{proof}

\begin{lemma}[Uniform offspring-law bounds and weak limits]
\label{lem:moments}
Under \condref{B}{(B)}, \condref{T}{(C)}, and \condref{D}{(D)}, there
are constants $m_*,m^*,p_*,v_*>0$ and an integer $R\ge0$ such that the
following three assertions hold.

\emph{(i) } For every $r\ge1$,
\begin{equation}\label{eq:single-moment-bounds}
 m_*\le m_r\le m^*,\qquad \Prob(Y_r\ge1)\ge p_*.
\end{equation}

\emph{(ii) } From every sequence $\{Y_{r_j}:j\ge1\}$ with
$r_j\uparrow\infty$, one can extract a further subsequence, still
denoted by $\{Y_{r_j}\}$, such that $Y_{r_j}\Rightarrow Y_*$, where
$Y_*$ is non-negative integer-valued with $\E Y_*^2\le K_2$,
\begin{equation}\label{eq:weighted-conv}
 \sum_{k\ge1}k^2
    |\Prob(Y_{r_j}=k)-\Prob(Y_*=k)|\longrightarrow0,
\end{equation}
and $\E Y_*\ge m_*$.

\emph{(iii) } For every $r>R$,
\begin{equation}\label{eq:eventual-laws}
 \Var(Y_r)\ge v_*,\qquad
 \gcd\{k\ge1:\Prob(Y_r=k)>0\}=1.
\end{equation}
\end{lemma}

\begin{proof}
Since $m_r=\mu_r/\mu_{r-1}$, condition \condref{B}{(B)} gives
$m_*=\underline\mu/\overline\mu$ and
$m^*=\overline\mu/\underline\mu$ as possible bounds. It also gives
\begin{equation}\label{Min1}
 \frac{\underline\mu}{\overline\mu}\le M_{i,n}
        \le\frac{\overline\mu}{\underline\mu}.
\end{equation}
By Cauchy--Schwarz and Lemma~\ref{lem:dini},
\[
 m_r^2=\bigl(\E[Y_r;Y_r\ge1]\bigr)^2
       \le\E Y_r^2\,\Prob(Y_r\ge1)
       \le K_2\Prob(Y_r\ge1).
\]
Taking $p_*=m_*^2/K_2>0$, we obtain (i).

We next prove (ii). By Lemma~\ref{lem:dini}, the bound
$\sup_r\E Y_r^2\le K_2$ implies that the family $\{Y_r:r\ge1\}$ is tight.
Given any sequence $r_j\uparrow\infty$, pass to a
weakly convergent further subsequence, still denoted by $Y_{r_j}$, with
$Y_{r_j}\Rightarrow Y_*$. Since $\Nzero$ is closed, $Y_*$ is also
supported on $\Nzero$ by \cite[Theorem~2.1]{Billingsley}. We now prove
\eqref{eq:weighted-conv}. For every fixed $K$,
\[
 \sum_{k=1}^K k^2\Prob(Y_*=k)
       =\lim_{j\to\infty}\sum_{k=1}^K k^2\Prob(Y_{r_j}=k)
       \le K_2.
\]
Letting $K\to\infty$, monotone convergence gives
\begin{equation}\label{Y*2}
 \E Y_*^2\le K_2.
\end{equation}
Given $\varepsilon>0$, Lemma~\ref{lem:dini} allows us to choose $K$ so
large that
\begin{equation}\label{Y*01}
 \sup_j\E[Y_{r_j}^2;Y_{r_j}>K]\le V(K)<\varepsilon
\end{equation}
and
\begin{equation}\label{Y*00}
 \E[Y_*^2;Y_*>K]<\varepsilon.
\end{equation}
For this fixed $K$,
\[
 \sum_{1\le k\le K}k^2
   \bigl|\Prob(Y_{r_j}=k)-\Prob(Y_*=k)\bigr|\longrightarrow0,
   \qquad j\to\infty.
\]
Splitting the full series at $k\le K$ and $k>K$, and controlling both
tails by \eqref{Y*00} and \eqref{Y*01}, we obtain
\[
 \limsup_{j\to\infty}\sum_{k\ge1}k^2
   \bigl|\Prob(Y_{r_j}=k)-\Prob(Y_*=k)\bigr|
 \le2\varepsilon.
\]
Since $\varepsilon>0$ is arbitrary, \eqref{eq:weighted-conv} follows.
Consequently,
\[
 |\E Y_{r_j}^2-\E Y_*^2|
 \le\sum_{k\ge1}k^2|\Prob(Y_{r_j}=k)-\Prob(Y_*=k)|
 \longrightarrow0,
\]
and therefore
\[
 |\E Y_{r_j}-\E Y_*|
 \le\sum_{k\ge1}k\,|\Prob(Y_{r_j}=k)-\Prob(Y_*=k)|
 \longrightarrow0,
\]
which yields $\E Y_*\ge m_*$.

We prove (iii) by contradiction. If the first assertion failed, we could
choose strictly increasing indices $r_j$ with $\Var(Y_{r_j})\to0$.
Tightness would give a weakly convergent further subsequence, still
denoted by $Y_{r_j}$, with limit $Y_*$. By \eqref{eq:weighted-conv},
$\Var(Y_*)=0$, contrary to \condref{D}{(D)}. Hence there are an integer
$R_1$ and $v_*>0$ such that $\Var(Y_r)\ge v_*$ for every $r>R_1$.

For the second assertion, suppose instead that there are strictly
increasing indices $r_j$ with $d_+(Y_{r_j})\ge2$. Passing to a weakly
convergent further subsequence gives $Y_{r_j}\Rightarrow Y_*$, and
\condref{D}{(D)} yields $d_+(Y_*)=1$. Since the gcd of a nonempty set of
positive integers is attained by a finite subset, there are positive
integers $k_1,\ldots,k_\ell$ such that
\[
 \Prob(Y_*=k_u)>0\quad(1\le u\le\ell),
 \qquad \gcd(k_1,\ldots,k_\ell)=1.
\]
For each fixed $k_u$, weak convergence gives
$\Prob(Y_{r_j}=k_u)\to\Prob(Y_*=k_u)$. Hence
$\Prob(Y_{r_j}=k_u)>0$ for all sufficiently large $j$. Thus
$d_+(Y_{r_j})$ divides $\gcd(k_1,\ldots,k_\ell)=1$, a contradiction.
Therefore, there is an integer $R_2$ such that
$\gcd\{k\ge1:\Prob(Y_r=k)>0\}=1$ for every $r>R_2$. Taking
$R=\max\{R_1,R_2\}$ proves \eqref{eq:eventual-laws}.
\end{proof}

Recall $a_{i,n}$ and $Q_{i,n}$ defined in \eqref{Qin} and \eqref{eq:a-def}.

\begin{lemma}[Uniform scale growth and survival bounds]\label{lem:survival}
Under \condref{B}{(B)}, \condref{T}{(C)}, and \condref{D}{(D)}, the
following two assertions hold.

\emph{(i) } There are positive constants $c,C$ such that for all $0\le i\le n$,
\begin{equation}\label{eq:a-growth}
 c(n-i+1)\le1+a_{i,n}\le C(n-i+1).
\end{equation}

\emph{(ii) }   $Q_{i,n}>0$ for all
$0\le i\le n$, and there are positive constants $c,C$ such that
\begin{equation}\label{eq:survival-bounds}
 \frac{c}{n-i+1}\le Q_{i,n}\le\frac{C}{n-i+1}.
\end{equation}
\end{lemma}

\begin{proof}
We first prove (i). Let
$\sigma_r=\Var(Y_r)/m_r^2$. Since
$\sigma_r=\nu_r+1/m_r-1$, we have
\begin{align}
 \frac{1+2a_{i,n}}{M_{i,n}}
 =\frac{\mu_i}{\mu_n}
       +\mu_i\sum_{r=i+1}^n\frac{\nu_r}{\mu_{r-1}}
 =1+\mu_i\sum_{r=i+1}^n\frac{\sigma_r}{\mu_{r-1}}.
 \label{eq:theta}
\end{align}
The last equality uses
$\mu_i\sum_{r=i+1}^n(1/\mu_{r-1}-1/\mu_r)=1-\mu_i/\mu_n$.
Together with \eqref{eq:eventual-laws} and \eqref{Min1} we obtain (i).

We turn to (ii). The Cauchy--Schwarz inequality and
\eqref{eq:process-moments} in Lemma~\ref{lem:factorial-moments} give
\[
 Q_{i,n}\ge\frac{(\E Z_{i,n})^2}{\E Z_{i,n}^2}
 =\frac{M_{i,n}}{1+2a_{i,n}}>0.
\]
Combining this with \eqref{Min1} and \eqref{eq:a-growth} gives the
lower bound.

 For the upper bound, write $h=n-i$. By (i), fix $c>0$ such that
$a_{i,n}\ge ch$ for all sufficiently large $h$, uniformly in $i$.
This constant is fixed independently of the truncation level $K$.
Lemma~\ref{lem:kersting2} gives, for every integer $K\ge2$,
\begin{equation}\label{eq:shape-real-lower}
 \varphi_r(s)\ge\frac{\nu_r}{2K}
                   -\frac{V(K)}{2K m_*^2},\qquad 0\le s\le1.
\end{equation}
By \cite[Lemma~5]{Kersting2020},
\begin{equation}\label{eq:reciprocal}
 \frac1{1-F_{i,n}(s)}
 =\frac1{M_{i,n}(1-s)}
   +\mu_i\sum_{r=i+1}^n
       \frac{\varphi_r(F_{r,n}(s))}{\mu_{r-1}},\qquad 0\le s<1.
\end{equation}
Taking $s=0$ and using \eqref{eq:shape-real-lower} gives
\begin{align*}
 \frac1{Q_{i,n}}
 &\ge\frac1{M_{i,n}}+\frac{a_{i,n}}{K M_{i,n}}
       -\frac{V(K)}{2K m_*^2}
          \mu_i\sum_{r=i+1}^n\frac1{\mu_{r-1}}\\
 &\ge\frac1{M_{i,n}}+
       \frac hK\left(c\frac{\underline\mu}{\overline\mu}
            -\frac{\overline\mu V(K)}{2\underline\mu m_*^2}\right),
\end{align*}
for all sufficiently large $h$. Here we used \eqref{Min1},
$a_{i,n}\ge ch$, and
$\mu_i\sum_{r=i+1}^n\mu_{r-1}^{-1}
\le(\overline\mu/\underline\mu)h$.
By \condref{T}{(C)}, $V(K)\to0$. We can therefore choose one integer
$K\ge2$, independently of $i,n$, such that
\[
 V(K)\le c m_*^2
              \left(\frac{\underline\mu}{\overline\mu}\right)^2.
\]
For this fixed $K$,
\[
 \frac1{Q_{i,n}}
 \ge\frac{c\underline\mu}{2K\overline\mu}h,
\]
which gives the upper bound in \eqref{eq:survival-bounds} for all
sufficiently large $h$. The remaining bounded values of $h$ are
covered by increasing $C$, since $Q_{i,n}\le1$. For $i=n$,
$Q_{n,n}=1$, and the conclusion holds as well.
\end{proof}

\subsection{Single-step shape functions and composition errors}

In this subsection, we estimate composition errors of the shape
functions. For an offspring law $Y$ with $\E Y>0$, $\E Y^2<\infty$,
and generating function $f$, write $m=\E Y$ and
$p_+=\Prob(Y\ge1)$. For $s\in[0,1)$, define
\begin{align}\label{eq:HG}
 H(s)& =\frac{1-f(s)}{1-s}=\sum_{j\ge0}\Prob(Y>j)s^j,
\nonumber\\
G(s)&=\frac{f(s)-1+m(1-s)}{(1-s)^2}
       =\sum_{j\ge0}\E[(Y-j-1)_+]s^j.
\end{align}
For $s=1$, set $H(1):=m$ and $G(1):=f''(1)/2$.

\begin{lemma}[The positive-support gcd criterion]
\label{lem:gcd-inversion}
Suppose that $H$ is defined by \eqref{eq:HG}. Then the following three
conditions are equivalent:

(i) $d_+(Y)=1$;

(ii) $H(s)\ne0$ for every $s\in\overline{\dd}$;

(iii) the reciprocal $1/H$ has absolutely summable coefficients.
\end{lemma}

\begin{proof}
(i)$\Rightarrow$(ii).\quad By definition, $H(1)=m\neq0$. For $s\ne1$,
$H(s)=(1-f(s))/(1-s)$. Since $m>0$, for $|s|<1$ we have
$|f(s)|\le f(|s|)<1$, and therefore $H(s)\ne0$.

It remains to consider $|s|=1$ with $s\ne1$. Write
$p_k=\Prob(Y=k)$. If $f(s)=1$, taking real parts gives
\begin{equation}\label{eq:real-part}
 \sum_{k\ge1}p_k(1-\Re s^k)=0.
\end{equation}
Since $|s|=1$ implies $\Re s^k\le1$, every term in
\eqref{eq:real-part} is nonnegative. Thus $\Re s^k=1$ for every $k$
with $p_k>0$, and hence $s^k=1$. Since $m>0$, at least one such $k$
exists, so $s$ is a root of unity; let $d\ge2$ be its order. Now
$s^k=1$ is equivalent to $d\mid k$, so $d$ divides every positive
support point. Therefore,
\[
 d\ \Big|\ \gcd\{k\ge1:p_k>0\}=d_+(Y),
\]
contradicting $d_+(Y)=1$. This proves (ii).

(ii)$\Rightarrow$(i).\quad If $d:=d_+(Y)>1$, let
$\zeta=\exp(2\pi\mathrm i/d)$. Every positive support point is a
multiple of $d$, so $f(\zeta)=1$. Since $\zeta\ne1$, it follows
that $H(\zeta)=0$.

(ii)$\Rightarrow$(iii).\quad This follows from
Lemma~\ref{lem:wiener}.

(iii)$\Rightarrow$(ii).\quad If $1/H$ has absolutely summable
coefficients, it is continuous on $\overline{\dd}$, and the
pointwise identity $H\cdot(1/H)=1$ shows that $H$ has no zero there.
\end{proof}

Define $H_r,G_r$ by \eqref{eq:HG} with $Y$ replaced by $Y_r$.
For $r>R$, Lemmas~\ref{lem:moments} and \ref{lem:gcd-inversion}
show that $1/H_r$ has absolutely summable coefficients. By
\eqref{varphi00} and \eqref{eq:HG},
\begin{equation}\label{varphir0}
 \varphi_r=m_r^{-1}G_rH_r^{-1},\qquad r>R,
\end{equation}
whose coefficient norm is finite by the convolution inequality.

\begin{lemma}[Relative compactness of inverses and shape functions]
\label{lem:shape-compact}
For every $r>R$, $H_r(s)\ne0$ for every $s\in\overline{\dd}$. The
families $\{H_r^{-1}:r>R\}$ and $\{\varphi_r:r>R\}$ are relatively
compact in the coefficient norm, where $R$ is from
Lemma~\ref{lem:moments}. In particular,
\begin{equation}\label{eq:shape-uniform}
 \sup_{r>R}\bigl(\norm{H_r^{-1}}+\norm{G_r}+\norm{\varphi_r}\bigr)<\infty.
\end{equation}
\end{lemma}

\begin{proof}
By \eqref{eq:HG},
\[
 \norm{H_r}=\sum_{j\ge0}\Prob(Y_r>j)=\E Y_r=m_r,\qquad
 \norm{G_r}=\sum_{j\ge0}\E[(Y_r-j-1)_+]=\frac{f_r''(1)}2,
\]
so $\norm{H_r},\norm{G_r}<\infty$. For every $r>R$,
Lemma~\ref{lem:moments}(iii) gives $d_+(Y_r)=1$.
Lemma~\ref{lem:gcd-inversion} therefore shows that $H_r$ is zero-free
on $\overline{\dd}$ and that $1/H_r$ has absolutely summable
coefficients.

Any bounded sequence of indices in $\{r:r>R\}$ has a constant
subsequence, along which both families trivially converge. Any
unbounded sequence of indices has an increasing subsequence
$\{r_j:j\ge1\}$. Take a further weakly convergent subsequence
$Y_{r_j}\Rightarrow Y_*$. Denote the corresponding functions of $Y_*$
by $H_*,G_*$, and write $p_k^{(r_j)}=\Prob(Y_{r_j}=k)$ and
$p_k^*=\Prob(Y_*=k)$. By \eqref{eq:weighted-conv},
\begin{align*}
 \norm{H_{r_j}-H_*}
 &\le\sum_{k\ge1}k|p_k^{(r_j)}-p_k^*|\longrightarrow0,\\
 \norm{G_{r_j}-G_*}
 &\le\frac12\sum_{k\ge2}k(k-1)|p_k^{(r_j)}-p_k^*|\longrightarrow0.
\end{align*}
Lemma~\ref{lem:moments}(ii) gives $0<m_*\le\E Y_*<\infty$, and
\condref{D}{(D)} gives $d_+(Y_*)=1$. Applying
Lemma~\ref{lem:gcd-inversion} shows that $H_*^{-1}$ has absolutely
summable coefficients. Let $D_j=H_*^{-1}(H_{r_j}-H_*)$. Then
$\norm{D_j}\to0$ and $H_{r_j}=H_*(1+D_j)$. For large $j$,
$\norm{D_j}<1$, so
\[
 \sum_{\ell\ge0}(-D_j)^\ell=(1+D_j)^{-1}
\]
and
\[
 H_{r_j}^{-1}=H_*^{-1}(1+D_j)^{-1}
 =H_*^{-1}\left(1+\sum_{\ell\ge1}(-D_j)^\ell\right).
\]
Consequently,
\[
 \norm{H_{r_j}^{-1}-H_*^{-1}}
 \le\norm{H_*^{-1}}\frac{\norm{D_j}}{1-\norm{D_j}}
 \longrightarrow0.
\]
Finally, the representation \eqref{varphir0} shows that
$\varphi_{r_j}$ converges in the coefficient norm, since the means
$m_{r_j}$ also converge by \eqref{eq:weighted-conv}. Thus every
sequence in either family has a convergent subsequence, proving
relative compactness and hence boundedness in the coefficient norm.
Together with $\norm{G_r}=f_r''(1)/2\le K_2/2$, this gives
\eqref{eq:shape-uniform}.
\end{proof}

\begin{lemma}[Second-moment error bound after composition]\label{lem:shape-modulus}
There exists $C>0$ such that, for every $r\ge1$ and $0\le t\le1$,
\begin{equation}\label{eq:shape-scalar-modulus}
 |\varphi_r(1-t)-\varphi_r(1)|\le C(\omega(t)+t),
\end{equation}
where $\omega(\cdot)$ is defined by \eqref{Omegat0}. Moreover, for any
probability generating function $F$, put $Q=1-F(0)$. Then, uniformly
over $r>R$,
\begin{equation}\label{eq:shape-modulus}
 \norm{\varphi_r\circ F-\varphi_r(1)}
 \le C\bigl(\omega(Q)+Q\bigr).
\end{equation}
\end{lemma}

\begin{proof}
Write $q=F(0)$. By \eqref{eq:composition-contraction},
\begin{equation}\label{eq:Fq-min}
 \norm{F^j-1}=2\bigl(1-q^j\bigr)\le2\min\{jQ,1\},\qquad j\ge0,
\end{equation}
where $q^0=1$. By \eqref{eq:HG}, the sum of the coefficients of $G_r$
is $G_r(1)=f_r''(1)/2$, so $G_r\circ F$ has finite coefficient norm by
\eqref{eq:composition-contraction}. For every integer $k\ge0$,
\[
 \sum_{j\ge0}(k-j-1)_+=\frac{k(k-1)}2,
\]
and $(k-j-1)_+>0$ implies $j\le k$. Therefore,
\begin{align}\label{eq:G-first}
 \norm{G_r\circ F-G_r(1)}
 &\le2\sum_{j\ge0}\E[(Y_r-j-1)_+]\min\{jQ,1\}\nonumber\\
 &\le\E\bigl[Y_r(Y_r-1)\min\{Y_rQ,1\}\bigr]
 \le\omega(Q).
\end{align}
Similarly, using $\sum_{j\ge0}j\Prob(Y_r>j)=f_r''(1)/2$,
\begin{align}\label{eq:Hr0}
 \norm{H_r\circ F-m_r}
 \le2Q\sum_{j\ge0}j\,\Prob(Y_r>j)
 =Qf_r''(1)\le K_2 Q.
\end{align}

For the real estimate, fix any $r\ge1$ and set $x=1-t\in[0,1]$.
Coefficient calculations similar to those above, together with
$|x^j-1|\le\min\{jt,1\}$, give
\[
 |G_r(x)-G_r(1)|\le\omega(t),\qquad
 |H_r(x)-m_r|\le K_2t.
\]
Since the coefficients of $H_r$ are nonnegative,
\eqref{eq:single-moment-bounds} yields
$H_r(x)\ge H_r(0)=\Prob(Y_r\ge1)\ge p_*$. Using
$\varphi_r(x)=G_r(x)/(m_rH_r(x))$ and
\eqref{eq:single-moment-bounds}, we obtain
\begin{align*}
 |\varphi_r(x)-\varphi_r(1)|
 &\le\frac{|G_r(x)-G_r(1)|}{m_rH_r(x)}
   +\frac{G_r(1)}{m_r}
         \left|\frac1{H_r(x)}-\frac1{m_r}\right|\\
 &\le\frac{\omega(t)}{m_*p_*}
   +\frac{K_2}{2m_*}\frac{|m_r-H_r(x)|}{m_*p_*}
 \le C(\omega(t)+t).
\end{align*}
This proves \eqref{eq:shape-scalar-modulus}.

For the rest of the proof, fix $r>R$. By
\eqref{eq:composition-product}, $H_r^{-1}\circ F$ is the pointwise
reciprocal of $H_r\circ F$, and $H_r^{-1}(1)=m_r^{-1}$. Hence
\[
 H_r^{-1}\circ F-H_r^{-1}(1)
 =\frac{m_r-H_r\circ F}{m_r}\,\bigl(H_r^{-1}\circ F\bigr).
\]
The composition contraction, \eqref{eq:shape-uniform}, \eqref{eq:Hr0},
and $m_r\ge m_*$ now yield
\begin{align}\label{eq:Hr1}
 \norm{H_r^{-1}\circ F-m_r^{-1}}
 \le m_r^{-1}\,\norm{H_r^{-1}}\cdot Qf_r''(1)\le C Q.
\end{align}

Finally, from $\varphi_r=m_r^{-1}G_rH_r^{-1}$ and
\eqref{eq:composition-product}, adding and subtracting
$m_r^{-1}G_r(1)(H_r^{-1}\circ F)$ gives the decomposition
\[
 \varphi_r\circ F-\varphi_r(1)
 =m_r^{-1}\bigl(G_r\circ F-G_r(1)\bigr)\bigl(H_r^{-1}\circ F\bigr)
 +m_r^{-1}G_r(1)\bigl(H_r^{-1}\circ F-m_r^{-1}\bigr).
\]
Taking norms and using $G_r(1)=\norm{G_r}=f_r''(1)/2$,
$\norm{H_r^{-1}\circ F}\le\norm{H_r^{-1}}$,
\eqref{eq:shape-uniform}, \eqref{eq:G-first}, and \eqref{eq:Hr1}, we
obtain \eqref{eq:shape-modulus} uniformly over $r>R$.
\end{proof}

By Lemma~\ref{lem:survival}, fix $C_0\ge1$ such that
\begin{equation}\label{Qin000}
 Q_{i,n}\le\min\{1,C_0/(n-i+1)\}.
\end{equation}
For integers $h\ge0$, define
\begin{equation}\label{eq:rhotilde}
 \epsilon_h:=\omega\!\left(\min\left\{1,\frac{C_0}{h+1}\right\}\right)
                 +\frac1{h+1},\qquad
 S_h:=\sum_{\ell=0}^{h-1}\epsilon_\ell,\qquad S_0=0.
\end{equation}
By Lemma~\ref{lem:dini}, as $h\to\infty$,
\begin{equation}\label{eq:eps-summability}
 \epsilon_h\longrightarrow0,\qquad
 \sum_{h\ge0}\frac{\epsilon_h}{h+1}<\infty,\qquad S_h=o(h).
\end{equation}
The first assertion follows from $\omega(t)\to0$ as $t\downarrow0$.
For the second, \eqref{eq:dini-sum} handles the modulus term, and
$\sum_h(h+1)^{-2}<\infty$ handles the added term. For the last
assertion, fix $\delta>0$ and choose $H$ such that
$\epsilon_\ell\le\delta$ for $\ell\ge H$. Then, for $h>H$,
\[
 \frac{S_h}{h}\le
 \frac1h\sum_{\ell=0}^{H-1}\epsilon_\ell+\delta.
\]
Letting $h\to\infty$ and then $\delta\downarrow0$ gives the desired
result.

\begin{lemma}[Linear fractional approximation]\label{lem:linear-fractional}
For every $0\le i\le n$, let
\[
 U_{i,n}(s)=1-F_{i,n}(s),\qquad
 L_{i,n}(s)=\frac{M_{i,n}(1-s)}{1+a_{i,n}(1-s)}, \qquad 0\le s<1.
\]
Define
\[
 \RR_{i,n}(s)=\mu_i\sum_{r=i+1}^n
   \frac{\varphi_r(F_{r,n}(s))-\varphi_r(1)}{\mu_{r-1}},\qquad 0\le s<1.
\]
(i) With $h=n-i$, we have
\begin{equation}\label{eq:UL-bound}
 \norm{U_{i,n}}+\norm{L_{i,n}}\le\frac{C}{h+1},
 \qquad |\RR_{i,n}(0)|\le CS_h.
\end{equation}
(ii) \begin{equation}\label{eq:survival-asymptotic}
 \lim_{N\to\infty}
 \sup_{\substack{0\le i<n\\n-i\ge N}}
 \left|\frac{Q_{i,n}a_{i,n}}{M_{i,n}}-1\right|=0.
\end{equation}
(iii)
For  $i\ge R$, the remainder $\RR_{i,n}$ has
absolutely summable coefficients, and
\begin{equation}\label{eq:reciprocal-remainder}
 \frac1{U_{i,n}(s)}=\frac1{L_{i,n}(s)}+\RR_{i,n}(s),
 \qquad |s|<1,
\end{equation}
with
\begin{equation}\label{eq:UL-error-tail}
 \norm{U_{i,n}-L_{i,n}}\le\frac{CS_h}{(h+1)^2},\qquad i\ge R.
\end{equation}
\end{lemma}

\begin{proof}
(i) Since $\varphi_r(1)=\nu_r/2$, by the definition of $a_{i,n}$, we have
\[
 \mu_i\sum_{r=i+1}^n\frac{\varphi_r(1)}{\mu_{r-1}}
  =\frac{\mu_i}{\mu_n}\frac{\mu_n}{2}
           \sum_{r=i+1}^n\frac{\nu_r}{\mu_{r-1}}
  =\frac{a_{i,n}}{M_{i,n}}.
\]
Also, for $0\le s<1$,
\[
 \frac1{L_{i,n}(s)}
  =\frac1{M_{i,n}(1-s)}+\frac{a_{i,n}}{M_{i,n}}.
\]
Recalling \eqref{eq:reciprocal}, we obtain
\begin{equation}\label{eq:real-reciprocal-remainder}
 \frac1{U_{i,n}(s)}=\frac1{L_{i,n}(s)}+\RR_{i,n}(s),
 \qquad 0\le s<1,
\end{equation}
for every starting point.

The norm estimates follow from the signs of the coefficients.
The constant coefficient of $U_{i,n}=1-F_{i,n}$ is $Q_{i,n}$;
all other coefficients are nonpositive, with absolute sum $Q_{i,n}$.
Thus
\[
 \norm{U_{i,n}}=2Q_{i,n}.
\]
For $a=a_{i,n}\ge0$, the geometric expansion is
\begin{equation}\label{ain0}
 \frac1{1+a(1-s)}=\frac1{1+a}
       \sum_{j\ge0}\left(\frac{a}{1+a}\right)^j s^j.
\end{equation}
If $a>0$, this series has radius $1+1/a>1$; if $a=0$, the
function equals $1$, and no expression $1/a$ is used.
Write $\gamma_j=a^j/(1+a)^{j+1}$, $j\ge0$, with
$\gamma_0=1$ and $\gamma_j=0$ for $j\ge1$ when $a=0$.
These coefficients are nonnegative, nonincreasing and sum to one.
With $\gamma_{-1}=0$, by the definition of $L_{i,n}(s)$ and \eqref{ain0} we have
\[
 L_{i,n}(s)=M_{i,n}(1-s)\sum_{j\ge0}\gamma_js^j
          =M_{i,n}\sum_{j\ge0}(\gamma_j-\gamma_{j-1})s^j.
\]
Its constant coefficient is positive and all later coefficients are
nonpositive. Then
\begin{align*}
 \norm{L_{i,n}}
 =M_{i,n}\gamma_0+
          M_{i,n}\sum_{j\ge1}(\gamma_{j-1}-\gamma_j)
 =2M_{i,n}\gamma_0=\frac{2M_{i,n}}{1+a_{i,n}}.
\end{align*}
Using \eqref{eq:survival-bounds}, \eqref{Min1} and
\eqref{eq:a-growth}, we get the first bound in \eqref{eq:UL-bound}.

By \eqref{eq:shape-scalar-modulus} of Lemma~\ref{lem:shape-modulus} and
$q_{r,n}=1-Q_{r,n}$,
\[
 |\varphi_r(q_{r,n})-\varphi_r(1)|
          \le C(\omega(Q_{r,n})+Q_{r,n}),\qquad r\ge1.
\]
By \eqref{Qin000}, \eqref{eq:rhotilde} and monotonicity of $\omega$,
\begin{align*}
 \omega(Q_{r,n})+Q_{r,n}
 \le\omega\!\left(\min\left\{1,\frac{C_0}{n-r+1}\right\}\right)
            +\frac{C_0}{n-r+1}\le(1+C_0)\epsilon_{n-r}.
\end{align*}
Since $\mu_i/\mu_{r-1}\le\overline\mu/\underline\mu$,
\begin{equation}\label{eq:R-scalar-bound}
 |\RR_{i,n}(0)|\le C\sum_{r=i+1}^n\epsilon_{n-r}
             =C\sum_{\ell=0}^{h-1}\epsilon_\ell=CS_h.
\end{equation}
This proves the second bound in \eqref{eq:UL-bound}.

(ii) Taking $s=0$ in \eqref{eq:real-reciprocal-remainder} yields
\[
\left|\frac{Q_{i,n}a_{i,n}}{M_{i,n}}-1\right|
=\left|\frac{1+M_{i,n}\RR_{i,n}(0)}{1+a_{i,n}+M_{i,n}\RR_{i,n}(0)}\right|.
\]
For all sufficiently large $h$, Lemma~\ref{lem:survival}(i) gives
$a_{i,n}\ge ch>0$. By \eqref{eq:eps-summability} and
\eqref{eq:R-scalar-bound}, as $h\to \infty$, we have
\[
 \left|\frac{1+M_{i,n}\RR_{i,n}(0)}{a_{i,n}}\right|
             \le C\frac{1+S_h}{h}\longrightarrow0
\]
 uniformly in $i<n$, which yields \eqref{eq:survival-asymptotic}.

(iii) For $r>i\ge R$,   $\varphi_r$ has finite coefficient norm by
Lemma~\ref{lem:shape-compact}, and so does $\varphi_r\circ F_{r,n}$
by \eqref{eq:composition-contraction}. The finite sum defining
$\RR_{i,n}$ therefore has absolutely summable coefficients.
The coefficient-norm estimate in Lemma~\ref{lem:shape-modulus} gives
\[
 \norm{\varphi_r\circ F_{r,n}-\varphi_r(1)}
                  \le C(\omega(Q_{r,n})+Q_{r,n})
                  \le C\epsilon_{n-r}.
\]
Summing with the same bounded weights yields
\begin{equation}\label{eq:R-bound}
 \norm{\RR_{i,n}}\le C\sum_{r=i+1}^n\epsilon_{n-r}=CS_h,
 \qquad i\ge R.
\end{equation}

For $|s|<1$, $|F_{i,n}(s)|\le F_{i,n}(|s|)<1$. Hence $1/U_{i,n}$ is
analytic in the open disc. The only zero of $L_{i,n}$ in the closed
disc is at $s=1$. Consequently
$1/L_{i,n}$ is also analytic in the open disc.
Both sides of \eqref{eq:real-reciprocal-remainder} are now analytic
there and agree on $(0,1)$, so they coincide on the open disc; this
proves
\eqref{eq:reciprocal-remainder}. Multiplication by $U_{i,n}L_{i,n}$
gives
\begin{equation}\label{eq:UL-identity}
 U_{i,n}-L_{i,n}=-U_{i,n}L_{i,n}\RR_{i,n}
\end{equation}
pointwise on $(0,1)$. Both sides have absolutely summable
coefficients, so this is an equality of power series and holds
pointwise on the whole closed disc. Taking norms and applying
\eqref{eq:UL-bound} and \eqref{eq:R-bound},
\[
 \norm{U_{i,n}-L_{i,n}}
  \le\norm{U_{i,n}}\norm{L_{i,n}}\norm{\RR_{i,n}}
  \le\frac{C}{h+1}\frac{C}{h+1}CS_h
  =\frac{CS_h}{(h+1)^2}.
\]
This proves \eqref{eq:UL-error-tail}. When $i=n$, the sum is empty,
$U_{n,n}=L_{n,n}=1-s$ and the error is zero.
\end{proof}

\subsection{Second-order Taylor expansion and the exact product recursion}

For $i\le n$ define the conditional generating function and the
correction factor
\begin{equation}\label{eq:correction-def}
 \Psi_{i,n}(s):=\frac{F_{i,n}(s)-q_{i,n}}{Q_{i,n}},\qquad
 \HH_{i,n}(s):=[1+a_{i,n}(1-s)]\Psi_{i,n}(s).
\end{equation}
Since $Q_{i,n}>0$, $\Psi_{i,n}$ is a probability generating function
concentrated on the positive integers, $\norm{\Psi_{i,n}}=1$,
$\norm{\HH_{i,n}}\le1+2a_{i,n}<\infty$, and $\HH_{i,n}(1)=1$. These definitions also cover
$i=n$, in which case $\Psi_{n,n}=\HH_{n,n}=s$.
The coefficients of $\HH_{i,n}$ need not be nonnegative; all norm and
coefficient-tail estimates for these factors use absolute values.
Since $\Psi_{r,n}(0)=0$, we define the one-step factor by an explicit
series rather than by division. Fix $1\le r\le n$, and write
$F=F_{r,n}$, $q=F(0)$ and $Q=1-q>0$. Put
\begin{equation}\label{eq:alpha-definition}
 d_{r,n}=\frac{Q_{r-1,n}}{Q_{r,n}}=
        \frac{1-f_r(q)}Q,
 \qquad
 \alpha_{r,n}(s)=\frac1{d_{r,n}}
       \sum_{k\ge1}\Prob(Y_r=k)
          \sum_{j=0}^{k-1}q^{k-1-j}F(s)^j.
\end{equation}
This series $\alpha_{r,n}(s)$ has finite coefficient norm: its coefficients
are nonnegative, and the sum of its coefficients before normalization
is
\[
 \sum_{k\ge1}\Prob(Y_r=k)\sum_{j=0}^{k-1}q^{k-1-j}
     =\frac{1-f_r(q)}{1-q}=d_{r,n}\le m_r<\infty.
\]
Thus $\norm{\alpha_{r,n}}=\alpha_{r,n}(1)=1$.
Factoring $F^k-q^k=(F-q)\sum_{j=0}^{k-1}q^{k-1-j}F^j$ gives
\begin{equation}\label{eq:psi-recursion}
 \Psi_{r-1,n}=\alpha_{r,n}\Psi_{r,n}.
\end{equation}
Finally, with $U_{r,n}=1-F_{r,n}$ and $\beta_r=f_r''(1)/m_r$, put
\begin{equation}\label{eq:alpha-expansion}
 b_{r,n}=\alpha_{r,n}-1+\frac{\beta_r}{2}U_{r,n}.
\end{equation}

\begin{lemma}[Expansion of the one-step conditional factor]
\label{lem:one-step}
For every $1\le r\le n$, the functions $\alpha_{r,n}$ and $b_{r,n}$
have absolutely summable coefficients. With $h=n-r$,
\begin{align}
 \norm{b_{r,n}}
 &\le C\bigl(Q_{r,n}^2+Q_{r,n}\omega(Q_{r,n})\bigr)
 \le C\left[\frac1{(h+1)^2}+\frac{\epsilon_h}{h+1}\right],
 \label{eq:b-bound}\\
 \norm{\alpha_{r,n}-1}&\le\frac{C}{h+1}.
 \label{eq:alpha-one-bound}
\end{align}
These estimates include bounded remaining times and initial generations.
\end{lemma}

\begin{proof}
Abbreviate $f=f_r$, $m=m_r$, $\gamma=f_r''(1)$,
$F=F_{r,n}$, $q=q_{r,n}$, $Q=Q_{r,n}$ and $d=d_{r,n}$.
Let $\Delta=F-q$. Its coefficients are nonnegative,
$\Delta(0)=0$, and $\norm\Delta=\Delta(1)=Q$.
Since $q<1$, each derivative $f^{(\ell)}(q)$ is finite. Define
\[
 T(s)=\sum_{\ell=2}^{\infty}
            \frac{f^{(\ell)}(q)}{\ell!}\Delta(s)^{\ell-2}.
\]
All coefficients in this series are nonnegative. By Taylor expansion
at $q$, with its value at one obtained by monotone convergence, and
then the integral remainder formula,
\begin{align*}
 \norm T=T(1)
 =\frac{f(1)-f(q)-Qf'(q)}{Q^2}
 =\frac1{Q^2}\int_q^1(1-t)f''(t)\,\mathrm dt
 \le\frac\gamma2.
\end{align*}
Hence $\norm T<\infty$.
For real $s\in[0,1]$, the expansion of $f(q+\Delta(s))$ has
nonnegative terms, so Tonelli's theorem justifies the
rearrangement
\begin{align}
 f(F(s))-f(q)
 =f'(q)\Delta(s)+\sum_{\ell=2}^{\infty}
                 \frac{f^{(\ell)}(q)}{\ell!}\Delta(s)^\ell
 =\Delta(s)\bigl[f'(q)+\Delta(s)T(s)\bigr].
 \label{eq:taylor-factor}
\end{align}
Since $\norm\Delta=Q$ and
$\sum_{\ell\ge1}f^{(\ell)}(q)Q^\ell/\ell!=f(1)-f(q)<\infty$, all
series here have absolutely summable coefficients. Summing
the identities
$F^k-q^k=\Delta\sum_{j=0}^{k-1}q^{k-1-j}F^j$ with weights
$\Prob(Y_r=k)$ gives
\[
f\circ F-f(q)=\Delta\,d\alpha_{r,n}.
\]
For $s\in(0,1)$ we have
$\Delta(s)>0$; comparison with \eqref{eq:taylor-factor} and
absolute summability of the coefficients therefore give
\begin{equation}\label{eq:alpha-T}
 d\alpha_{r,n}=f'(q)+\Delta T.
\end{equation}
Since $d>0$, $\alpha_{r,n}$ has absolutely summable coefficients;
by \eqref{eq:alpha-expansion}, so does $b_{r,n}$.

To quantify the Taylor error, we compare $T$ with $\gamma/2$.
The constant coefficient of $T$ is $f''(q)/2$ and every other
coefficient is nonnegative. The bound $T(1)\le\gamma/2$ gives
\begin{equation}\label{tdelta}
 \norm{T-\gamma/2}
 =\frac{\gamma-f''(q)}2+T(1)-\frac{f''(q)}2
 \le\gamma-f''(q).
\end{equation}
For $j\ge0$, $1-q^j\le\min\{jQ,1\}$. Therefore
\begin{align}
 0\le\gamma-f''(q)
 &=\sum_{k\ge2}k(k-1)\Prob(Y_r=k)(1-q^{k-2})\nonumber\\
 &\le\sum_{k\ge2}k^2\Prob(Y_r=k)\min\{kQ,1\}
 \le\omega(Q).
 \label{eq:second-tail}
\end{align}
Combining the last two estimates,
\begin{equation}\label{eq:T-modulus}
 \norm{T-\gamma/2}\le\omega(Q).
\end{equation}
It remains to account for the normalization by $d$. Evaluating
\eqref{eq:alpha-T} at $s=1$ and using $\alpha_{r,n}(1)=1$ and
$\Delta(1)=Q$ give
\begin{equation}\label{eq:d-def}
 d=f'(q)+QT(1).
\end{equation}
We need a positive bound for $d$ on all intervals. Directly from
\eqref{eq:alpha-definition},
\[
 d=\sum_{k\ge1}\Prob(Y_r=k)(1+q+\cdots+q^{k-1})
       \ge\Prob(Y_r\ge1)\ge p_*.
\]
Also, $d=Q^{-1}\int_q^1f'(t)\,\mathrm dt$, so
\[
 m-d=\frac1Q\int_q^1\bigl(m-f'(t)\bigr)\,\mathrm dt
     =\frac1Q\int_q^1(t-q)f''(t)\,\mathrm dt.
\]
Since $0\le f''(t)\le\gamma\le K_2$, these estimates give
\begin{equation}\label{eq:d-bounds}
 d\ge p_*,\qquad 0\le m-d\le\frac\gamma2Q\le CQ.
\end{equation}
By \eqref{eq:alpha-T}, \eqref{eq:d-def} and using
$\Delta=Q-U_{r,n}$ give
\begin{align*}
 d(\alpha_{r,n}-1)
 =\Delta T-QT(1) =-\frac\gamma2U_{r,n}+\Delta(T-\gamma/2)
       -Q\bigl(T(1)-\gamma/2\bigr).
\end{align*}
Using \eqref{eq:alpha-expansion} and $\beta_r=\gamma/m$, we obtain
\begin{equation}\label{eq:b-exact}
 b_{r,n}=
    \left(\frac\gamma{2m}-\frac\gamma{2d}\right)U_{r,n}
       +\frac{\Delta(T-\gamma/2)
              -Q\bigl(T(1)-\gamma/2\bigr)}d.
\end{equation}
For its first term, \eqref{eq:d-bounds} gives
\[
 \left|\frac\gamma{2m}-\frac\gamma{2d}\right|
       =\frac{\gamma|d-m|}{2md}\le CQ,
 \qquad \norm{U_{r,n}}=2Q.
\]
Its norm is therefore at most $CQ^2$. For the second term, since
$|T(1)-\gamma/2|\le\norm{T-\gamma/2}\le\omega(Q)$,
\[
 \frac{\norm{\Delta(T-\gamma/2)
             -Q\bigl(T(1)-\gamma/2\bigr)}}d
 \le\frac{Q\norm{T-\gamma/2}+Q|T(1)-\gamma/2|}{p_*}
 \le CQ\omega(Q).
\]
Adding both terms proves
\begin{equation}\label{eq:b-prefinal}
 \norm{b_{r,n}}\le C(Q^2+Q\omega(Q)).
\end{equation}
The survival bound and the definition of $\epsilon_h$ give
\[
 Q^2\le\frac{C}{(h+1)^2},\qquad
 Q\omega(Q)\le\frac{C_0}{h+1}
       \omega\!\left(\min\left\{1,\frac{C_0}{h+1}\right\}\right)
       \le\frac{C_0\epsilon_h}{h+1}.
\]
This proves \eqref{eq:b-bound}. Finally,
$\alpha_{r,n}-1=-(\beta_r/2)U_{r,n}+b_{r,n}$, while
$\beta_r\le K_2/m_*$ and $\omega(Q)\le K_2$. Hence
\[
 \norm{\alpha_{r,n}-1}
       \le CQ+C(Q^2+Q\omega(Q))\le CQ\le\frac{C}{h+1},
\]
proving \eqref{eq:alpha-one-bound}.
\end{proof}

For all $1\le r\le n$, define the product error by
\begin{equation}\label{eq:e-definition}
 e_{r,n}=\left(1+\frac{\beta_r}{2}L_{r,n}\right)\alpha_{r,n}-1.
\end{equation}

\begin{lemma}[Uniform summability of the product error]\label{lem:product}
The functions $e_{r,n}$ have absolutely summable coefficients and
satisfy the exact recursion
\begin{equation}\label{eq:H-recursion}
 \HH_{r-1,n}=\HH_{r,n}(1+e_{r,n}),\qquad 1\le r\le n,
\end{equation}
where $\HH_{r-1,n}$ is defined by \eqref{eq:correction-def}. For $h=n-r\ge0$, let
\begin{equation}\label{eq:e-bound}
 \kappa_h=C\left[
   \frac1{(h+1)^2}+\frac{\epsilon_h}{h+1}
                  +\frac{S_h}{(h+1)^2}\right].
\end{equation}
With $C$ fixed sufficiently large and independently of $r,n$,
\begin{equation}\label{eq:e-two-ranges}
 \norm{e_{r,n}}\le
 \begin{cases}
   \kappa_{n-r},&r>R,\\
   C/(n-r+1),&1\le r\le R.
 \end{cases}
\end{equation}
Moreover,
\begin{equation}\label{eq:kappa-tail}
 \sum_{h\ge0}\kappa_h<\infty,\qquad
 \tau_N:=\sum_{h\ge N}\kappa_h\longrightarrow0.
\end{equation}
In particular, the total norm error is uniformly bounded over
finite time intervals:
\begin{equation}\label{eq:last}
 \sup_{0\le i<n}\sum_{r=i+1}^n\norm{e_{r,n}}
 \le\sum_{h\ge0}\kappa_h+CR<\infty.
\end{equation}
\end{lemma}

\begin{proof}
By \eqref{eq:a-def},
\[
 a_{r-1,n}-a_{r,n}
 =\frac{\mu_n}{2}\sum_{j=r}^n\frac{\nu_j}{\mu_{j-1}}
   -\frac{\mu_n}{2}\sum_{j=r+1}^n\frac{\nu_j}{\mu_{j-1}}
 =\frac{\mu_n\nu_r}{2\mu_{r-1}}.
\]
Since $\nu_r=\beta_r/m_r=\beta_r\mu_{r-1}/\mu_r$,
\begin{equation}\label{eq:a-increment}
 a_{r-1,n}-a_{r,n}=\frac{\beta_rM_{r,n}}2.
\end{equation}
Use \eqref{eq:correction-def} and \eqref{eq:psi-recursion} to write
\begin{align*}
 \HH_{r-1,n}
 &=\bigl[1+a_{r-1,n}(1-s)\bigr]\Psi_{r-1,n}\\
 &=\bigl[1+a_{r,n}(1-s)
           +(a_{r-1,n}-a_{r,n})(1-s)\bigr]
                \Psi_{r,n}\alpha_{r,n}\\
 &=\HH_{r,n}
       \left(1+\frac{\beta_r}{2}
            \frac{M_{r,n}(1-s)}{1+a_{r,n}(1-s)}\right)
                \alpha_{r,n}\\
 &=\HH_{r,n}\left(1+\frac{\beta_r}{2}L_{r,n}\right)\alpha_{r,n}.
\end{align*}
The last factor is $1+e_{r,n}$ by definition, proving
\eqref{eq:H-recursion}. Substitution of
$\alpha_{r,n}=1-(\beta_r/2)U_{r,n}+b_{r,n}$ gives
\begin{equation}\label{eq:e-exact}
 e_{r,n}=
 \frac{\beta_r}{2}(L_{r,n}-U_{r,n})+b_{r,n}
 -\frac{\beta_r^2}{4}L_{r,n}U_{r,n}
 +\frac{\beta_r}{2}L_{r,n}b_{r,n}.
\end{equation}
All these terms have absolutely summable coefficients for every
generation.

For $r>R$, applying \eqref{eq:UL-error-tail} and
$0\le\beta_r=f_r''(1)/m_r\le K_2/m_*$, we have
\[
 \norm{\tfrac{\beta_r}{2}(L_{r,n}-U_{r,n})}
                 \le C\frac{S_h}{(h+1)^2}.
\]
Meanwhile, Lemma~\ref{lem:one-step} gives
\[
 \norm{b_{r,n}}\le C\left[
                    \frac1{(h+1)^2}+\frac{\epsilon_h}{h+1}\right].
\]
For the third term of \eqref{eq:e-exact}, the global norm bound \eqref{eq:UL-bound} and
the convolution inequality \eqref{eq:convolution} give
\[
 \norm{\tfrac{\beta_r^2}{4}L_{r,n}U_{r,n}}
          \le C\norm{L_{r,n}}\norm{U_{r,n}}
          \le\frac{C}{(h+1)^2}.
\]
Finally, combining \eqref{eq:b-bound} with the bound for $L_{r,n}$,
\begin{align*}
 \norm{\tfrac{\beta_r}{2}L_{r,n}b_{r,n}}
 \le\frac{C}{h+1}
                \left[\frac1{(h+1)^2}+\frac{\epsilon_h}{h+1}\right]
 \le C\left[\frac1{(h+1)^2}+\frac{\epsilon_h}{h+1}\right].
\end{align*}
Adding these estimates proves the first case of
\eqref{eq:e-two-ranges}.

For  $r\le R$, we estimate the
error directly. By \eqref{eq:e-definition} and
$\norm{\alpha_{r,n}}=1$, together with  \eqref{eq:alpha-one-bound} and \eqref{eq:UL-bound}, we get
\begin{align*}
 \norm{e_{r,n}}
 \le\norm{\alpha_{r,n}-1}
         +\frac{\beta_r}{2}\norm{L_{r,n}}\norm{\alpha_{r,n}}
 \le\frac{C}{h+1},
\end{align*}
and \eqref{eq:e-two-ranges} is proved.

To complete the proof, we recall \eqref{eq:e-bound}.
The first term of \eqref{eq:e-bound} is summable. The second is
summable by \eqref{eq:eps-summability}. For the third, Tonelli's
theorem and $\epsilon_\ell\ge0$ give
\begin{align*}
 \sum_{h\ge1}\frac{S_h}{(h+1)^2}
 =\sum_{h\ge1}\frac1{(h+1)^2}\sum_{\ell=0}^{h-1}\epsilon_\ell
 =\sum_{\ell\ge0}\epsilon_\ell
              \sum_{h\ge\ell+1}\frac1{(h+1)^2}
 \le\sum_{\ell\ge0}\epsilon_\ell
              \int_{\ell+1}^\infty\frac{\mathrm dx}{x^2}
  =\sum_{\ell\ge0}\frac{\epsilon_\ell}{\ell+1}<\infty.
\end{align*}
The term $h=0$ has $S_0=0$ and is finite. Hence
\eqref{eq:kappa-tail} holds.   Finally, \eqref{eq:e-two-ranges}
and the fact that every interval contains at most $R$ initial
indices lead to \eqref{eq:last}.
\end{proof}

\subsection{Relative compactness of the family of correction factors}

Recall from \eqref{eq:correction-def} that
$\HH_{i,n}(s)=[1+a_{i,n}(1-s)]\Psi_{i,n}(s)$, where
$\Psi_{i,n}(s)=(F_{i,n}(s)-q_{i,n})/Q_{i,n}$.

\begin{proposition}[Relative compactness of the correction factors]
\label{prop:correction-compact}
The family
\[
 \mathscr H=\{\HH_{i,n}:0\le i\le n<\infty\}
\]
is relatively compact in the coefficient norm. In particular, there is
$C_H<\infty$ such that
\begin{equation}\label{eq:H-uniform}
 \sup_{0\le i\le n<\infty}\norm{\HH_{i,n}}\le C_H.
\end{equation}
Moreover, writing $\HH_{i,n}(s)=\sum_{j\ge0}c_j^{(i,n)}s^j$, one has
\begin{equation}\label{eq:H-tails}
 T_J:=\sup_{0\le i\le n<\infty}
        \sum_{j>J}|c_j^{(i,n)}|\longrightarrow0,
        \qquad J\to\infty.
\end{equation}
\end{proposition}

\begin{proof}
Since $\HH_{n,n}(s)=s$, iteration of \eqref{eq:H-recursion} gives
\begin{equation}\label{eq:H-product}
 \HH_{i,n}(s)=s\prod_{r=i+1}^{n}(1+e_{r,n}(s)),
 \qquad 0\le i\le n,
\end{equation}
with the convention $\prod_{r=n+1}^{n}=1$. Using $\norm{s}=1$ and
Lemma~\ref{lem:product}, we obtain, for $0\le i\le n$,
\begin{align}
 \norm{\HH_{i,n}}
 \le\prod_{r=i+1}^{n}(1+\norm{e_{r,n}})
 \le\exp\left\{\sum_{r=i+1}^{n}\norm{e_{r,n}}\right\}
 \le\exp\left\{\sum_{h\ge0}\kappa_h+CR\right\}=:C_H<\infty.
 \label{eq:product-exp}
\end{align}
This proves \eqref{eq:H-uniform}.

By Lemma~\ref{lem:dini},
\[
 \sup_{r\ge1}\sum_{k>J}\Prob(Y_r=k)\le\frac{K_2}{(J+1)^2};
\]
the coefficients at indices $0,\ldots,J$ lie in the compact set
$[0,1]^{J+1}$. Consequently $\{f_r:r\ge1\}$ is
totally bounded in the coefficient norm. Denote its compact closure by
$\mathscr K$. This closure consists of probability generating
functions, since coefficient
nonnegativity and the sum of coefficients equal to one are preserved
under $\ell^1$ convergence.

We also verify continuity of composition. Suppose $g_j\to g$
and $F_j\to F$ in the coefficient norm, with $F_j,F$ probability
generating functions. For every integer $k\ge1$,
\[
 F_j^k-F^k=(F_j-F)\sum_{\ell=0}^{k-1}F_j^\ell F^{k-1-\ell},
 \qquad \norm{F_j^k-F^k}\le k\norm{F_j-F}.
\]
Truncate $g$ at degree $J$ and use the composition contraction:
\begin{align*}
 \norm{g_j\circ F_j-g\circ F}
 &\le\norm{g_j-g}+\norm{g\circ F_j-g\circ F}\\
 &\le\norm{g_j-g}
      +\sum_{k=0}^J |g_k|k\norm{F_j-F}
      +2\sum_{k>J}|g_k|.
\end{align*}
For fixed $J$, the first two terms tend to zero as $j\to\infty$.
Then let $J\to\infty$, proving continuity of composition.
Repeated composition is therefore continuous on $\mathscr K^\ell$
for every fixed $\ell\ge1$. For each fixed $N\ge1$, the functions
$F_{i,n}$ with $n-i\le N$ belong to the compact set
\[
 \{s\}\ \cup\
 \bigcup_{\ell=1}^{N}
 \{g_1\circ\cdots\circ g_\ell:g_1,\ldots,g_\ell\in\mathscr K\}.
\]
Each set in the union is a continuous image of $\mathscr K^\ell$;
the function $s$ accounts for intervals of length zero. Consequently,
for each fixed $N$, the family
$\{F_{i,n}:0\le i\le n,\ n-i\le N\}$ is relatively compact.

On these intervals, $Q_{i,n}\ge c/(N+1)$ and $a_{i,n}\le C(N+1)$.
Let $(i_j,n_j)_{j\ge1}$ be any sequence of intervals with
$0\le i_j\le n_j$ and $n_j-i_j\le N$. By the relative compactness
just established, there are a subsequence $(j_k)$ and a limit $F$
such that $F_{i_{j_k},n_{j_k}}\to F$ in the coefficient norm.
Passing to a further subsequence, still denoted by $j_k$, we may
also assume that $a_{i_{j_k},n_{j_k}}\to a$ for some real number $a$.
Evaluation at zero is continuous in the coefficient norm, so
$q_{i_{j_k},n_{j_k}}\to F(0)$ and
$Q_{i_{j_k},n_{j_k}}\to1-F(0)$, with the latter bounded away from
zero. The operations in \eqref{eq:correction-def} are continuous
under these convergences, so $\HH_{i_{j_k},n_{j_k}}$ converges in
the coefficient norm. This proves relative compactness of
$\{\HH_{i,n}:0\le i\le n,\ n-i\le N\}$ for each fixed $N$.

We now extend this compactness result to arbitrary interval lengths
by truncating the product. For any integer $N\ge1$ and $n-i\ge N$,
iterate \eqref{eq:H-recursion} from $i$ to $n-N$ to obtain
\[
 \HH_{i,n}-\HH_{n-N,n}
 =\HH_{n-N,n}
          \left[\prod_{r=i+1}^{n-N}(1+e_{r,n})-1\right].
\]
Every remaining time in this product satisfies $n-r\ge N$. By
\eqref{eq:e-two-ranges}, the factors with $r>R$ have total norm error
at most $\tau_N$, while the at most $R$ initial factors have total norm
error at most $CR/(N+1)$. Expanding a finite product over its nonempty
subsets and applying the convolution inequality \eqref{eq:convolution}
gives
\[
 \norm{\prod_r(1+e_{r,n})-1}
       \le\prod_r(1+\norm{e_{r,n}})-1
       \le\exp\left\{\sum_r\norm{e_{r,n}}\right\}-1.
\]
By \eqref{eq:H-uniform}, $\norm{\HH_{n-N,n}}\le C_H$, with $C_H$
independent of $N$. Therefore,
\begin{equation}\label{eq:H-approx}
 \norm{\HH_{i,n}-\HH_{n-N,n}}
 \le C_H\left[\exp\left\{\tau_N+\frac{CR}{N+1}\right\}-1\right]
       \longrightarrow0\qquad(N\to\infty)
\end{equation}
uniformly over $0\le i\le n-N$.

Given $\varepsilon>0$, choose $N$ so that the right-hand side is less
than $\varepsilon$. Relative compactness of
$\{\HH_{i,n}:0\le i\le n,\ n-i\le N\}$ provides a finite
$\varepsilon$-net. By \eqref{eq:H-approx}, every $\HH_{i,n}$ from an
interval $0\le i\le n-N$ is within $\varepsilon$ of one from an
interval of length $N$, and hence within $2\varepsilon$ of this finite
net. Thus $\mathscr H$ is totally bounded. Since power series with
finite coefficient norm form a complete space, the closure of
$\mathscr H$ is compact.

Finally, let $\varepsilon>0$. Since $\mathscr H$ is totally bounded,
it has a finite $\varepsilon$-net $\{G^{(1)},\ldots,G^{(m)}\}$.
Write $G^{(\ell)}(s)=\sum_{j\ge0}g_j^{(\ell)}s^j$. Because there are
only finitely many centres, we can choose $J$ so large that
\[
 \max_{1\le\ell\le m}\sum_{j>J}|g_j^{(\ell)}|<\varepsilon.
\]
For any $\HH_{i,n}\in\mathscr H$, choose $\ell$ such that
$\norm{\HH_{i,n}-G^{(\ell)}}<\varepsilon$. Then
\[
 \sum_{j>J}|c_j^{(i,n)}|
 \le\norm{\HH_{i,n}-G^{(\ell)}}
       +\sum_{j>J}|g_j^{(\ell)}|
 <2\varepsilon.
\]
Taking the supremum over $0\le i\le n<\infty$ gives
$T_J\le2\varepsilon$. Since $\varepsilon$ is arbitrary, this proves
\eqref{eq:H-tails}.
\end{proof}

\section{Proof of Theorem~\ref{thm:main}}

\begin{proof}[Proof of Theorem~\ref{thm:main}]
Fix $i<n$ and abbreviate $a=a_{i,n}$, $Q=Q_{i,n}$, $M=M_{i,n}$, and
$h=n-i$. Write the coefficients of $\HH_{i,n}$ as
$c_j\equiv c_j^{(i,n)}$. By \eqref{eq:correction-def},
\begin{equation}\label{eq:geometric-factor}
 \Psi_{i,n}(s)=\frac{\HH_{i,n}(s)}{1+a(1-s)}.
\end{equation}
The coefficients of $(1+a(1-s))^{-1}$ are nonnegative, sum to one,
and are bounded by $1/(1+a)$. Write $[s^k]F$ for the coefficient of
$s^k$ in $F$. The convolution inequality \eqref{eq:convolution},
Lemmas~\ref{lem:moments} and \ref{lem:survival}, and
\eqref{eq:H-uniform} give
\[
 \sup_{k\ge1}p_{i,n}(k)
 =Q\sup_{k\ge1}\left|[s^k]
       \frac{\HH_{i,n}(s)}{1+a(1-s)}\right|
 \le\frac{Q}{1+a}\norm{\HH_{i,n}}
 \le\frac{C}{(h+1)^2}.
\]
This proves \eqref{eq:main-upper}.

We now prove the local asymptotic \eqref{eq:main-local}. By
Lemma~\ref{lem:survival}(i), choose a fixed $H_1$ such that
$a_{i,n}\ge ch>0$ whenever $h=n-i\ge H_1$, uniformly in $i$. For the
remainder of the proof, put $\theta=a/(1+a)\in(0,1)$.

Expanding
$(1+a(1-s))^{-1}=(1+a)^{-1}\sum_{m\ge0}\theta^ms^m$ gives
\begin{equation}\label{eq:exact-coefficient}
 p_{i,n}(k)=\frac{Q}{1+a}\theta^k
                \sum_{j=0}^{k}c_j\theta^{-j},\qquad k\ge1.
\end{equation}
Fix $L_0>0$ and restrict to $k\le L_0h$. Since $a\ge ch$, for
$0\le j\le k$,
\[
 1\le\theta^{-j}=(1+1/a)^j
       \le\mathrm e^{k/a}\le\mathrm e^{L_0/c}.
\]
For fixed $J$ and $k\ge J$, since
$\sum_{j\ge0}c_j=\HH_{i,n}(1)=1$,
\begin{align}\label{eq:coeff-split}
 \sum_{j=0}^{k}c_j\theta^{-j}-1
  =\sum_{j=0}^{J}c_j\bigl(\theta^{-j}-1\bigr)
  \quad+\sum_{j=J+1}^{k}c_j\bigl(\theta^{-j}-1\bigr)
       -\sum_{j>k}c_j .
\end{align}
The absolute value of the first sum is at most
$C_H[(1+1/a)^J-1]$. The remaining terms are controlled by the uniform
coefficient tail $T_J$ from \eqref{eq:H-tails}, giving
\begin{equation}\label{eq:coefficient-convergence}
 \left|\sum_{j=0}^{k}c_j\theta^{-j}-1\right|
 \le C_H\bigl[(1+1/a)^J-1\bigr]
       +(\mathrm e^{L_0/c}+1)T_J.
\end{equation}
Given $\varepsilon>0$, use \eqref{eq:H-tails} to choose a fixed $J$
such that $(\mathrm e^{L_0/c}+1)T_J<\varepsilon/2$. Then choose
$N\ge\max\{H_1,J\}$ large enough that
\[
 C_H\bigl[(1+1/(cN))^J-1\bigr]<\varepsilon/2.
\]
For every interval with $h\ge N$ and every integer
$N\le k\le L_0h$, we have $k\ge J$ and $a\ge cN$. Thus the right-hand
side of \eqref{eq:coefficient-convergence} is less than $\varepsilon$.
This proves that the coefficient sum tends to one uniformly in the
joint range of \eqref{eq:main-local}.

The survival asymptotic \eqref{eq:survival-asymptotic} and
$a/(1+a)\to1$ also give
\[
 \frac{Q a^2}{M(1+a)}
    =\frac{Qa}{M}\frac{a}{1+a}\longrightarrow1
\]
uniformly in $i$ as $h=n-i\to\infty$. Moreover, since
$0\le x-\log(1+x)\le x^2/2$ for $x\ge0$, the restrictions
$k\le L_0h$ and $a\ge ch$ yield
\[
 0\le\frac{k}{a}+k\log\theta
   =k\left[\frac1a-\log(1+1/a)\right]
   \le\frac{k}{2a^2}\le\frac{L_0}{2c^2h}\longrightarrow0.
\]
By \eqref{eq:exact-coefficient}, the normalized probability is the
product of the three quantities just estimated:
\[
 \frac{a^2}{M}\mathrm e^{k/a}p_{i,n}(k)
 =\frac{Qa^2}{M(1+a)}
   \exp\!\left\{k\left[\frac1a-\log(1+1/a)\right]\right\}
   \sum_{j=0}^k c_j\theta^{-j}.
\]
By above discussion, each factor tends to $1$ uniformly as $N\to\infty$ over
$h\ge N$ and $N\le k\le L_0h$. This proves \eqref{eq:main-local}.
The assertion for a process starting at time zero follows by taking
$i=0$. The proof is completed.
\end{proof}

{\bf Acknowledgement.}  The second author is supported by National Natural Science Foundation of China  (Grant No. 12671167).

\end{document}